\documentclass[12pt]{article}
\usepackage{theorem}
\usepackage{amssymb}
\usepackage{amsmath}

\newtheorem{thm}{Theorem}[section]
\newtheorem{prop}[thm]{Proposition}
\newtheorem{cor}[thm]{Corollary}
\newtheorem{lm}[thm]{Lemma}

{\theorembodyfont{\rm}

\newtheorem{defn}[thm]{Definition}

\newtheorem{rem}[thm]{Remark}

}

\def\ra{\rightarrow}

\def\C{{\Bbb C}}

\def\P{{\Bbb P}}

\def\C{{\Bbb C}}

\title{Density of rational points on \\
Enriques surfaces}

\author{F. A.  Bogomolov\\
\small  Courant Institute of Mathematical Sciences, N.Y.U. \\
\small 251 Mercer str. \\
\small New York, (NY) 10012, U.S.A.\\
\small e-mail: bogomolo@cims.nyu.edu\\
\small  and\\
Yu. Tschinkel\\
\small Dept. of Mathematics, U.I.C.\\
\small 851 South Morgan str.\\
\small Chicago, (IL) 60607-7045,  U.S.A.  \\
\small e-mail: yuri@math.uic.edu
}

\begin{document}

\date{}

\maketitle

\thispagestyle{empty}

\pagebreak

\section{Introduction}

Let $X$ be an irreducible algebraic 
variety defined over a number field $K$.
We will say that rational points on $X$ are potentially dense
if there exists a finite extension $K'$ of $K$ such 
that the set of $K'$-rational points is Zariski--dense in $X$.
It seems reasonable to ask whether or
not rational points are potentially dense if
neither $X$ nor its unramified coverings admit a morphism 
onto a variety of general type. This question has an easy answer 
if the dimension of $X$ equals 1. Very little is known
in higher dimensions. 
Clearly, potential density holds for unirational varieties
and for abelian varieties. 
In \cite{harris-tschi-98} and \cite{bogomolov-tschi-98}
it is proved that rational points are potentially dense on
all Fano threefolds with the possible exception of
double covers of $\P^3$ ramified in smooth surfaces of degree $6$.

\

In this paper we will study the question of density
under the additional assumption that $X$ admits 
the structure of an elliptic fibration $\varphi :  X\ra B$
over some irreducible normal base $B$ of dimension $\ge 1$. As one
of the applications we
prove  that rational points on Enriques surfaces 
are potentially dense. 
An alternative idea to prove density
would have been to use the group of automorphisms ${\rm Aut}(X)$. 
Notice, however, that there exist Enriques surfaces with a
finite group ${\rm Aut}(X)$ (cf. \cite{barth-peters}).

\

{\small
{\bf Acknowledgements.} 
The second author is very grateful to 
Barry Mazur and Joe Harris for 
their ideas and suggestions.
The paper was written while both
authors were enjoying the hospitality of 
the Max-Planck-Institute in Bonn.
The first author was partially supported by the NSF.
The second author was partially supported by the NSA. 
We would like to thank the referee for detailed
comments helping us to improve the exposition.
}

\section{Elliptic fibrations - Generalities}

Let $K$ be a number field and
$X$ a smooth projective algebraic 
surface defined over $K$.
We say that $X$ admits a structure of an elliptic fibration 
if there exists a regular map
$\varphi :  X\ra B$ onto a smooth (irreducible) 
curve $B$ whose fibers are connected curves such that the generic
fiber is a smooth curve of genus  $1$. 
We denote by $X_b$ the
fiber over $b\in B$.
We will  say that $X$ admits a structure of a 
Jacobian elliptic fibration
if there exists a zero section $e : B\ra X$. 
It is well known (cf. \cite{barth-peters-vdv-84} 
or \cite{isk-shaf-96}) that 
to every elliptic fibration 
$\varphi : {\cal E}\ra B$ one can associate a
Jacobian elliptic fibration 
$\varphi_{\cal J}: {\cal J}={\cal J}({\cal E})\ra B$ (over the
same groundfield), which over
the generic point of the base $B$ is given by classes of divisors
of degree zero in the fiber. The zero section corresponds to the 
trivial class. There is a fiberwise action of 
${\cal J}$ on ${\cal E}$,
more precisely, a rational map
$$
\psi :  {\cal J}\times_B {\cal E}\ra {\cal E}
$$
which is regular in non-singular points of the fibers of 
${\cal J}$ and ${\cal E}$ and which induces a 
transitive action of ${\cal J}_b$ on ${\cal E}_b$ 
(for smooth fibers).

\subsection{Multisections}

Let $\varphi :  {\cal E}\ra B$ be an elliptic fibration.
A multisection $i : {\cal M}\hookrightarrow {\cal E}$ is 
an irreducible subvariety of ${\cal E}$ 
such that the map $\varphi\circ i : {\cal M}\ra B$ 
is finite and surjective. 
We will denote by $d=d({\cal M})$ the degree of
this projection.

\begin{defn}\label{defn:order-m}
Let $\varphi :  {\cal E}\ra B$ be an elliptic fibration.
A multisection ${\cal M}$ is said to be of {\em order $m$} 
if $m$ is the smallest positive
integer such that for any $b\in B$ and any pair of points 
$p_b,p_b'\in {\cal M} \cap {\cal E}_b$ the image of the zero-cycle
$p_b-p_b'$  in ${\cal J}_b$ is torsion of order $m$. 
\end{defn}

Note that {\em any} section of an elliptic fibration is
a multisection of order one.

Let $\Phi_m\subset {\cal J}$ be the subvariety of
of $m$-torsion points of ${\cal J}$
(all points $p_b$ such that  $m\cdot p_b=0$ 
in the group of rational points of the corresponding fiber 
${\cal J}_b$) which are not contained in 
the zero section $e(B)$. 

\begin{lm}
If ${\cal M}\subset \psi(\Phi_m \times_B {\cal M})$ 
then ${\cal M}$ is a multisection of order $m$.
\end{lm}

{\em Proof.} Tautology. 

\begin{defn}\label{dfn:non-torsion-section}
Let $\varphi  : {\cal E}\ra B$ 
be an elliptic fibration. 
A {\em saliently ramified multisection} 
of ${\cal E}$ is a multisection ${\cal M}$ 
which is ramified in a point $p_b$ which lies in 
a smooth (elliptic curve) fiber ${\cal E}_b$ of 
${\cal E}$.
\end{defn}

\begin{prop} \label{prop:local}
Let ${\cal M}$ be a saliently ramified 
multisection of ${\cal E}$. 
Then ${\cal M}$ is not of order $m$ for any $m>0$. 
\end{prop}

{\em Proof.} 
By assumption,
nearby fibers of ${\cal E}_b$ are smooth. Consider
an embedding $K\hookrightarrow \C$.    
We can find a sequence $(b_i)\in B(\C)$ of points
in the base converging to $b$ (in complex topology) 
and pairs of distinct points 
$p_i,p_i'$ in the fibers ${\cal E}_{b_i}(\C)$ 
which converge to $p$. 
If we assume that the cycle $p_i-p_i'$ is torsion of order $m$ 
(for some $m\ge 2$) in the Jacobian we obtain a contradiction 
since it converges to $0=e(b)$.

\begin{defn}\label{defn:tau} Let ${\cal M}$ be a multisection
of degree $d$ of the elliptic fibration $\varphi : {\cal E}\ra B$ 
and $b\in B$. We denote by 
${\rm Tr}_{\cal M}(b)$ the zero-cycle ${\cal M}\cap {\cal E}_b$.
Define the class map 
$$
\tau_{\cal M}\, :\,  {\cal E}_B \ra {\cal J}_B
$$
by the following rule:
$$
\tau_{\cal M}(p)=[d\cdot p -   {\rm Tr}_{\cal M}(\varphi(p))]
$$
for $p\in {\cal E}$. 
\end{defn}

\begin{lm}\label{lm:tau}
Suppose that the multisection ${\cal M}$ is of degree $d$ and 
not of order $d'$ with $d'|d$. 
Then the map 
$$
\tau_{\cal M}\,:\, {\cal M}\ra \tau_{\cal M}({\cal M})
$$
is a birational isomorphism. 
\end{lm}

{\em Proof.} Consider all points $p,p'\in {\cal E}$ such 
that $\tau_{\cal M}(p)=\tau_{\cal M}(p')$. Then the cycle
$p-p'$ is of order $d'|d$ in ${\cal J}_{b}$. On the other hand,
for a given $d'$ pairs of distinct points differing by torsion
of order $d'$ (more precisely, by a translation by $\Phi_{d'}$)
constitute a Weil divisor in ${\cal M}$.  Therefore,
the map $\tau_{\cal M}$ is a birational isomorphism.

\begin{prop}\label{prop:properties}
The map $\tau_{\cal M}$ is 
regular outside of singular fibers
of ${\cal E}\ra B$. 
\end{prop}

{\em Proof.} Evident.

\begin{cor} \label{cor:tau-of-M}
Let ${\cal M}$ be a saliently ramified 
multisection of ${\cal E}$.
Then $\tau_{\cal M}({\cal M})$  is a saliently 
ramified multisection of
${\cal J}$.
\end{cor}

{\em Proof.} 
The map $\tau_{\cal M}$ is unramified on the set
of non-singular fibers.  By lemma \ref{lm:tau}, 
it is a birational isomorphism. 

\

\begin{rem}\label{cor:tau-of-M-2}
Let ${\cal M}\subset {\cal E}$ be a multisection which is not
of order $m$ for any $m$. Then $\tau_{\cal M}({\cal M})$ is not
of order $m$ for any $m$. 
\end{rem}

From now on we will   
restict to the case of the base $B=\P^1$. Merel's theorem 
implies:

\begin{prop}\cite{bogomolov-tschi-98}\label{theo:non-torsion}
Let $ \varphi_{\cal J}  :  {\cal J}\ra \P^1$ be a Jacobian elliptic
fibration defined over $K$ 
with a saliently ramified multisection ${\cal M}$.
Then for all but finitely many 
$b\in \varphi_{\cal J}({\cal M}(K))\subset \P^1(K)$ 
the fibers $J_{b}$ have infinitely many rational points. 
\end{prop}

\begin{thm}\label{thm:ell}
Let $\varphi :  {\cal E}\ra \P^1$ be an elliptic fibration defined
over $K$. Assume that there exists 
a rational or elliptic saliently ramified multisection
${\cal M}$ defined over $K$. 
Then rational points on ${\cal E} $ are potentially 
dense. 
\end{thm}

{\em Proof.} Extending the groundfield, we can assume 
that ${\cal M}$ is defined over $K$, 
that $K$-rational points on ${\cal M}$ are Zariski-dense
and that it intersects a smooth fiber ${\cal E}_{\varphi(p)}$ 
with local intersection multiplicity $\ge 2$ at a 
$K$-rational point $p$. 
The image $\varphi({\cal M}(K))$ is Zariski-dense in $\P^1$. 
Every fiber ${\cal E}_b$ for $b\in \varphi({\cal M}(K))$ contains
at least one point $p_b\in {\cal M}(K)$. 
Therefore, it suffices to show that for almost all 
$b\in \varphi({\cal M}(K))$ the fiber ${\cal J}_b$ has a 
$K$-rational point of infinite order. Then we use the action of
the ${\cal J}$ on ${\cal E}$ to translate $p_b$. 

By assumption and \ref{cor:tau-of-M} the curve 
$\tau_{\cal M}({\cal M})$ is a saliently ramified 
multisection of ${\cal J}\ra \P^1$. 
The point $\tau_{\cal M}(p_b)\in \tau_{\cal M}({\cal M})$ is a
$K$-rational point of the fiber ${\cal J}_{\varphi(p_b)}$. 
Moreover, it is contained in the 
saliently ramified multisection $\tau_{\cal M}({\cal M})$. 
By theorem \ref{theo:non-torsion}, for all but
finitely many fibers ${\cal J}_{\varphi(p_b)}$ 
the point $\tau_{\cal M}(p_b)$ is a non-torsion point in the group 
${\cal J}_{\varphi(p_b)}(K)$. This concludes the proof.

\begin{rem}
An alternative argument avoiding Merel's theorem 
goes as follows: We can find
a base change $ \beta : B'\ra B$ with the following
properties: $\beta $ is \'etale at $b$ (corresponding to the
smooth fiber of ${\cal E}\ra B$ where ${\cal M}$ is ramified),
${\cal M}$ pulls back to a section ${\cal M}'$
of ${\cal E}'\ra B'$ and ${\cal E}'$ acquires a zero section $e'$
(which is different from ${\cal M}'$).
Then ${\cal M}'$ must be of infinite order in the Mordell--Weil 
group of ${\cal E}'$. A specialization argument 
(cf. \cite{silverman})
implies that rational points are (potentially) dense on ${\cal E}'$.
\end{rem}

\begin{prop}\label{prop:two-ell}
Let $X$ be a smooth algebraic surface defined over a number
field $K$ and admitting two different elliptic fibrations
$\varphi_1 : {\cal E}_1\ra B_1$ and 
$\varphi_2 : {\cal E}_2\ra B_2$ with $B_1\simeq B_2\simeq \P^1$.
Then rational points on $X$ are potentially dense. 
\end{prop}

{\em Proof.} If there is at least one ${\cal E}_{b_2}$ 
which is not of order $m$ 
for any $m$ we are done by \ref{cor:tau-of-M-2}. 
In view of \ref{cor:tau-of-M}, 
it suffices to consider the
case when for all $b_2\in B_2$ all multiple intersection
points of the multisections 
${\cal E}_{b_2}$ with fibers ${\cal E}_{b_1}$ (for $b_1\in B_1$) 
are contained in the singular fibers of the fibration ${\cal E}_1\ra B_1$. 
If a generic ${\cal E}_{b_2}$ is of some fixed order $m_0$ we see that
(a cover of) $X$ is dominated (birationally) by 
a product of two elliptic curves.

\section{Enriques surfaces}

We start with  a brief summary of 
the structure theory of Enriques surfaces 
(cf. \cite{barth-peters-vdv-84}, pp. 274-275):
Either $X$ has two distinct elliptic fibrations  over $\P^1$
(this case is called {\em non-special} in 
\cite{barth-peters-vdv-84})
or it admits an elliptic fibration with a $2$-section 
which is a $(-2)$-curve. 
In the second case (called {\em special})
the associated K3 double cover $\pi : Y\ra X$ 
admits the structure of an elliptic
fibration $Y\ra \P^1$ with {\em two} (non-intersecting) sections
(which could differ by torsion). 
The surface $Y$ is the minimal resolution of a double covering
of a quadratic cone $Q\subset \P^3$ (given by
$z_0z_1=z_2^2$ in standard coordinates in $\P^3=(z_0,z_1,z_2,z_3)$)
ramified in an intersection of this cone with a quartic 
hypersurface in $\P^3$ not passing through the vertex of this cone. 
The ramification curve  $R$ is reduced, of degree $8$
and has at most simple singularities. 
The fibers of the elliptic fibration correspond to the generators
of the cone $Q$ and the two sections 
are mapped to the vertex of $Q$ 
(cf. \cite{barth-peters-vdv-84}, p. 278).

\begin{thm} 
Let $X$ be an Enriques surface over $K$.
Then rational points  on $X$ are potentially dense.
\end{thm}

{\em Proof.}
In the non-special case we apply \ref{prop:two-ell}
and we are done.

Now let us consider the special case. 
Recall that the elliptic fibration $Y\ra \P^1$ has
two non-intersecting sections $e_1,e_2$. Their difference
is torsion in the Picard group of the generic point of 
$Y\ra \P^1$ if there exists an integer $m>0$
such that  $m(e_1-e_2)$ can be represented as a sum,
with integer coefficients, of components of 
singular fibers of $Y\ra \P^1$. If all fibers of
$Y\ra \P^1$ are irreducible then the difference $e_1-e_2$
is not torsion and rational points are potentially dense.
Otherwise, we have to consider subcases of the special case.

After blowing up the vertex of the cone
we can realize $Q$ birationally as 
$\P_{\P^1}({\cal O}(2)\oplus {\cal O})$
(see \cite{bogomolov-tschi-98}).
Taking sections of $\P_{\P^1}({\cal O}(2))$
(a 3-dimensional linear space) we obtain a family of (conics)
$\P^1_s\simeq \P^1\subset Q$ where each $\P^1_s$ intersects
the ramification curve $R$ in 8 
points (counted with multiplicities). 
The double cover $D_s$  of $\P^1_s$ is  
a multisection of the elliptic fibration $Y\ra \P^1$. We want
to find an elliptic saliently ramified multisection among the $D_s$.
Then we apply \ref{thm:ell}. 

Denote by $R^0\subset R$ the Zariski open 
subset of points where the curve
$R$ is smooth and where it is not 
tangent to the generators of the cone $Q$. 
For every point $P\in R^0$ we consider the affine line 
$L_P$ of sections $\P^1_s$ which are tangent to $R$ at $P$.
These sections cover the whole cone $Q$,
except the line joining $P$ and the vertex of $Q$. 
Consider the subset of sections of $\P_{\P^1}({\cal O}(2))$
which have at least two distinct points of 
local intersection multiplicity
$\ge 2$ with $R$. If this subset covers birationally the cone $Q$,
then we obtain a 1-dimensional family of elliptic curves
which covers $Y$ and which is generically transversal 
to the elliptic fibration $Y\ra \P^1$. Hence we can
apply \ref{prop:two-ell}. 

In particular, 
if the curve $R$ has a singular double point $r\in R$, 
then for any $P\in R^0$ we can find a 
tangent section $\P^1_s$ which 
passes through  $r$. 
The family of such sections covers $Q$ (birationally), provided that
$R$ has a component which is not a section. 
Therefore, we have only to consider the case
when $R$ consists of 4 distinct components which are sections.
Suppose that $r$ lies on the intersection of
two components. Then the tangents to some point of a third
component passing through $r$ cover $Q$ birationally and the
argument above applies. 

The only remaining case is the case of a non-singular curve $R$.
For any $P\in R$ there exists a section 
$s\in L_P$ which has local 
intersection multiplicity $\ge 2$ with $R$ 
at some other point. Indeed, we have 
a natural tangent correspondence 
${\rm TC}_P\subset L_P\times R$
defined by
$$
{\rm TC}_P = \{ (s,r)\,|\, s\in L_P, r\in \P^1_s\cap R\}.
$$
Since $R$ is irreducible this correspondence is irreducible. 
Since $R$ is not a rational curve this correspondence 
has ramification points over $L_P$. Ramification points correspond 
to sections having at least two intersection points of 
local multiplicity $\ge 2$ with $R$. 
Hence, the family of sections with 
this property (as $P$ moves over $R^0$)
covers (birationally) the cone $Q$ and defines a 
transversal elliptic fibration on $Y\ra \P^1$.


\begin{thebibliography}{1}

\bibitem{barth-peters}
W.~Barth and C.~Peters,
\newblock Automorphisms of {E}nriques surfaces,
\newblock {\em Invent. Math.}, 73(3):383--411, 1983.

\bibitem{barth-peters-vdv-84}
W.~Barth, C.~Peters, and A.~van~de Ven,
\newblock {\em Compact complex surfaces}, volume~4 of {\em Ergebnisse der
  Mathematik und ihrer Grenzgebiete (3)},
\newblock Springer-Verlag, Berlin, 1984.




\bibitem{bogomolov-tschi-98} 
F. A. Bogomolov and Yu. Tschinkel,
\newblock On the density of rational points on elliptic fibrations,
\newblock Preprint 1998.

\bibitem{harris-tschi-98}
J. Harris and Yu. Tschinkel,
\newblock Rational points on quartics,
\newblock Preprint 1998.

\bibitem{isk-shaf-96}
V.~A. Iskovskikh and I.~R. Shafarevich.
\newblock Algebraic surfaces,
\newblock In {\em Algebraic geometry, II}, volume~35 of {\em Encyclopaedia
  Math. Sci.}, pages 127--262. Springer, Berlin, 1996.


\bibitem{merel}
L. Merel,
\newblock 
Bornes pour la torsion des courbes elliptiques sur les corps de
nombres,
\newblock {\em Invent. Math.}, 124(1-3):437--449, 1996.

\bibitem{silverman} 
J. Silverman,
\newblock Heights and the specialization map for families of abelian
varieties,
\newblock {\em J. Reine Angew. Math.}, {342}:197--211, 1983.


\end{thebibliography}
\end{document}